\input amstex \loadmsbm
\documentstyle{amsppt}
\magnification=1200 \nopagenumbers \baselineskip=13pt
\parskip=4pt plus 1pt
\topskip=35pt
\hsize 15.5truecm
\vsize 22.5truecm

\font\ver=cmcsc10

\font\aut=cmr12

 \font\rms=cmr8 \font\rm=cmr9 \font\bfs=cmbx8 \font\bfm=cmbx9 \font\bfg=cmbx10 \font\its=cmti8
\font\itm=cmti9

 \font\tit=cmss17  \font\pizarra=msbm10

%%%%%%%%%%%%%%%%%%%%%%%%%%%%%%%%%%%%%%%

\def\headimpar{\hfil\ver The Leavitt path algebra of a graph  \hfil{\rm\folio}}
\def\headpar{{\rm \folio}\hfil\ver Abrams and Aranda \hfil}
\headline={\ifnum\pageno=1\hfil
\else\ifodd\pageno\headimpar\else\headpar\fi\fi}
\footline={\ifnum\pageno=1\hfil\tenrm{\rm \folio}\hfil\else\hfil\fi}

%%%%%%%%%%%%%  automatic numbering device
%%%%%%%%%%%%%%%%%%%%%%%%%%%%%%%%%%%%%%%
\newcount\parracount
\newcount\recount
\parracount=0
\recount=0
%%%%%%%%%%%%%%%%%%%%%%%%%%%%%%%%%%%%%

\def\re#1#2{\advance\recount by 1
\edef#1{\number\parracount .\number\recount}
\def\nnnn{#2}
\ifx\nnnn\empty\vskip 3 pt plus 1 pt{\bf #1. }\else
\vskip 3 pt plus 1 pt{\bf #1.}{ \ver #2. }\it\fi}
\def\endre{\par\rm}
\def\proof{{\bf Proof:} }
\def\endproof{~\vrule height 1.2ex width 1.2ex depth -0.06ex \par}
\def\parrafo#1{\advance\parracount by 1\recount=0
\vfil\vskip 10 pt plus 10 pt\noindent{\bf \number\parracount. #1}\par\nobreak\vskip 3 pt  plus 3 pt\nobreak}
%%%%%%%%%%%%%%%%%%%%%%%%%%%%%%%%%%%%%%%

%%%%%%%%%%%%%%%%%%%  Definitions
\def\s{\sigma}
\def\t{\tau}

%%%%%%%%%%%%%%%%%%%  References
%%%%%%%%%%%%%%%%%%%%%%%%%%%%%%%%
%\def\Leavitt{1}
%\def\Raeburn{2}
%\def\AGGP{3}

\def\AGGP{1}
\def\BPRS{2}
\def\C{3}
\def\CK{4}
\def\KPR{5}
\def\Lone{6}
\def\Ltwo{7}
\def\R{8}
\def\RS{9}

\endre

%%%%%%%%%%%%%%%%%%%%%%%%%  NOCIONES PRELIMINARES
%%%%%%%%%%%%%%%%%%%%%%%%%%%%%%%

\vskip 20 pt \centerline{\tit The Leavitt path algebra of a graph}

\vskip 14 pt \centerline{\aut Gene Abrams $^{a,}$\plainfootnote{\rms $^*$}{\rms Corresponding author.} \aut and Gonzalo
Aranda Pino $^b$} \vskip 9 pt \centerline{\its $^a$Department of Mathematics, University of Colorado, Colorado Springs,
CO 80933, U.S.A.} \vskip -1 pt \centerline{\its $^b$Departamento de \'Algebra, Geometr\'\i a y Topolog\'\i a,
Universidad de M\'alaga, 29071 M\'alaga, Spain} \plainfootnote{}{{\its E-mail addresses:} {\rms abrams\@math.uccs.edu
(G. Abrams), gonzalo\@agt.cie.uma.es (G. Aranda).}}

%\plainfootnote{}{{\its 2000 Mathematics Subject Classification: 16S10, 16S36, 16W50}}

\vskip 15 pt

{\narrower\noindent\baselineskip=5 pt {\bfs Abstract}

\rms For any row-finite graph $E$ and any field $K$ we construct the {\its Leavitt path algebra} $L(E)$ having
coefficients in $K$. When $K$ is the field of complex numbers, then $L(E)$ is the algebraic analog of the Cuntz Krieger
algebra $C^*(E)$ described in \cite{\R}.  The matrix rings $M_n(K)$ and the Leavitt algebras $L(1,n)$ appear as
algebras of the form $L(E)$ for various graphs $E$.  In our main result, we give necessary and sufficient conditions on
$E$ which imply that $L(E)$ is simple.}

\medskip {{\its Keywords:} {\rms path algebra; Leavitt algebra; Cuntz Krieger C*-algebra}}

\bigskip \noindent{\bf Introduction}\par\nobreak\vskip 3 pt  plus 3 pt\nobreak

Throughout this article $K$ will denote an arbitrary field. In his seminal paper \cite{\Lone}, Leavitt describes a
class of $K$-algebras (nowadays denoted by $L(m,n)$) which are universal with respect to an isomorphism property
between finite rank free modules. In \cite{\Ltwo}, Leavitt goes on to show that the algebras of the form $L(1,n)$ are
simple. More than a decade later, Cuntz \cite{\C} constructed and investigated the C*-algebras ${\Cal O}_n$ (nowadays
called the Cuntz algebras), showing, among other things, that each ${\Cal O}_n$ is (algebraically) simple.  When $K$ is
the field ${\Bbb C}$ of complex numbers, then ${\Cal O}_n$ can be viewed as the completion, in an appropriate norm, of
$L(1,n)$. Soon after the appearance of \cite{\C}, Cuntz and Krieger \cite{\CK} described the significantly more general
notion of the C*-algebra of a (finite) matrix $A$, denoted ${\Cal O}_A$.  Among this class of C*-algebras one can find,
for any finite graph $E$, the Cuntz-Krieger algebra $C^*(E)$, defined originally in \cite{\KPR}.  These C*-algebras, as
well as those arising from various infinite graphs, have been the subject of much investigation (see e.g. \cite{\R},
\cite{\RS}, and \cite{\BPRS}). Recently, the `algebraic analogs' of the C*-algebras  ${\Cal O}_A$ have been presented
in \cite{\AGGP}; these are denoted by ${\Cal {CK}}_A(K)$. By restricting attention to a specific set of allowable
matrices, the simplicity of the algebra ${\Cal {CK}}_A(K)$ for some subset of these allowable matrices has been
determined (although the condition for simplicity is not explicitly given in terms of the matrix $A$).

The goal of this article is to `complete the algebraic picture'. Specifically, we give the definition of the {\it
Leavitt path algebra} $L(E)$ corresponding to any row-finite graph $E$ and field $K$.   When $E$ is finite without
sources and sinks, then $L(E)$ can be realized as an algebra of the form ${\Cal {CK}}_A(K)$ for some matrix $A$.
Analogous to the relationship that exists between $L(1,n)$ and ${\Cal O}_n$, $L(E)$ has the property that when $K={\Bbb
C}$, then $C^*(E)$ can be viewed as the completion, in an appropriate norm, of $L(E)$ \cite{\R, Proposition 1.20}.

In our main result, Theorem (3.11), we give necessary and sufficient conditions on the row-finite graph $E$ which imply
that $L(E)$ is simple.   These results extend those presented in \cite{\AGGP}, in that:  they apply also to some
important algebras which are explicitly not considered in \cite{\AGGP}; they apply also to algebras which arise from
infinite matrices; and they provide necessary conditions on $E$ for the simplicity of $L(E)$. The statement of Theorem
(3.11) parallels a similar theorem for C*-algebras of the form $C^*(E)$ given in \cite{\R, Theorem 4.9 and subsequent
remarks}. However, the techniques utilized here are significantly different than those used in the analytic setting.

We begin by establishing some notational conventions. A {\bfm (directed) graph} $E=(E^0,E^1,r,s)$ consists of two
countable sets  $E^0,E^1$ and functions $r,s:E^1 \to E^0$. The elements of $E^0$ are called {\bfm vertices} and the
elements of $E^1$  {\bfm edges}. For each edge $e$, $s(e)$ is the {\bfm source} of $e$ and $r(e)$ is the {\bfm range}
of $e$. If $s(e)=v$  and $r(e)=w$, then we also say that $v$ {\bfm emits} $e$ and that $w$ {\bfm receives} $e$, or that
$e$ {\bfm points to} $w$.

A vertex which does not receive any edges is called a {\bfm source}. A vertex which emits no edges is called a {\bfm
sink}. A graph is called {\bfm row-finite} if $s^{-1}(v)$ is a finite set for each vertex $v$. In this paper, {\itm we
will only be concerned with row-finite graphs}. Of course, under this hypothesis, the edge set of $E$, $E^1$, is finite
if its set of vertices, $E^0$, is finite. Thus, we will say a graph $E$ is {\bfm finite} if $E^0$ is a finite set. A
{\bfm path} $\mu$ in a graph $E$ is a sequence of edges $\mu=\mu_1 \dots \mu_n$ such that $r(\mu_i)=s(\mu_{i+1})$ for
$i=1,\dots,n-1$. In such a case, $s(\mu):=s(\mu_1)$ is the source of $\mu$ and $r(\mu):=r(\mu_n)$ is the range of
$\mu$. If $s(\mu)=r(\mu)$ and $s(\mu_i)\neq s(\mu_j)$ for every $i\neq j$, then $\mu$  is a called a {\bfm cycle}.
\bigskip

\parrafo{Leavitt path algebras}

\medskip

In this section we define the algebraic structures under investigation. We begin by reminding the reader of the
construction of the standard path algebra of a graph.

\re{\pathalgebra}{}{\bf Definition.} Let $K$ be a field and $E$ be a graph.  The {\bfm path $K$-algebra over $E$} is
defined as the free $K$-algebra $K[E^0\cup E^1]$ with the relations: \item{(1)} $v_iv_j=\delta_{ij}v_i$ for every
$v_i,v_j\in E^0$. \item{(2)} $e_i=e_ir(e_i)=s(e_i)e_i$ for every $e_i\in E^1$.

This algebra is denoted by $A(E)$.
\endre

\re{\graphalgebra}{}{\bf Definition.} Given a graph $E$  we define the  {\bfm extended graph of} $E$ as the new graph
$\widehat{E}=(E^0,E^1\cup (E^1)^*,r',s')$ where $(E^1)^*=\{e_i^*:e_i\in  E^1\}$ and the functions $r'$ and $s'$ are
defined as $$r'|_{E^1}=r,\ s'|_{E^1}=s,\ r'(e_i^*)=s(e_i) \hbox{ and }  s'(e_i^*)=r(e_i).$$ \endre

\re{\graphalgebra}{}{\bf Definition.} Let $K$ be a field and $E$ be a  row-finite graph. The {\bfm Leavitt path algebra
of $E$ with coefficients in $K$} is defined as the path algebra over the  extended graph $\widehat{E}$, with relations:
\item{(CK1)} $e_i^*e_j=\delta_{ij}r(e_j)$ for every $e_j\in E^1$ and  $e_i^*\in (E^1)^*$. \item{(CK2)}
$v_i=\sum_{\{e_j\in E^1:s(e_j)=v_i\}}e_je_j^*$ for every  $v_i\in E^0$ which is not a sink.

This algebra is denoted by $L_K(E)$ (or more commonly simply by $L(E)$). \endre

The conditions CK1 and CK2 are called the {\bfm Cuntz-Krieger relations}. In  particular condition CK2 is the {\bfm
Cuntz-Krieger relation at $v_i$}. If $v_i$ is a sink, we do not have a CK2  relation at $v_i$. Note that the condition
of row-finiteness is needed in order to define the equation CK2.

\re{\examples}{}{\bf Examples.} Many well-known algebras are of the form $L(E)$ for some graph $E$:

\item{(i)} Matrix algebras $M_n(K)$: Consider the graph $E$ defined by $E^0=\{v_1,\dots,v_n\}$,
$E^1=\{e_1,\dots,e_{n-1}\}$ and $s(e_i)=v_i$ and $r(e_i)=v_{i+1}$ for $i=1,\dots,n-1$. Then $M_n(K)\cong L(E)$, via the
map $v_i \mapsto e(i,i), e_i \mapsto e(i,i+1)$, and $e_i^* \mapsto e(i+1,i)$ (where $e(i,j)$ denotes the standard
$(i,j)$-matrix unit in $M_n(K)$).
\item{(ii)} Laurent polynomial algebras $K[x,x^{-1}]$: Consider the graph $E$ defined by $E^0= \{*\}$, $E^1=\{x\}$. Then
clearly $K[x,x^{-1}]\cong L(E)$.
\item{(iii)} Leavitt algebras $A = L(1,n)$ for $n \geq 2$ investigated in \cite{\Ltwo}: Consider the graph $E$ defined by
$E^0=\{*\}$, $E^1=\{y_1,\dots,y_n\}$. Then $L(1,n)\cong L(E)$.
\endre

\re{\monomial}{\bf Lemma} Every monomial in $L(E)$ is of the following form. \item{(a)} $k_iv_i$ with $k_i\in K$ and
$v_i\in E^0$, or \item{(b)} $ke_{i_1}\dots e_{i_\s}e_{j_1}^*\dots e_{j_\t}^*$ where $k\in K;\  \s,\ \t\geq 0,\
\s+\t>0,\ e_{i_s}\in E^1$ and $e_{j_t}^*\in (E^1)^*$ for $0 \leq s \leq \sigma, 0 \leq t \leq \tau$. \endre

\proof  The proof is almost identical to the proof of \cite{\R, Corollary 1.15} (a straightforward induction argument
on the length of the monomial $kx_1\dots x_n$ with $x_i\in E^0\cup E^1 \cup (E^1)^*$), and so is omitted.
\endproof

\re{\unital}{\bf Lemma} If $E^0$ is finite then $L(E)$ is a unital  $K$-algebra. If $E^0$ is infinite, then $L(E)$ is
an algebra with local units (specifically, the set generated by finite sums of distinct elements of $E^0$). \endre

\proof First assume that $E^0$ is finite: we will show that $\sum_{i=1}^n  v_i$ is the unit element of the algebra.
First we compute $(\sum_{i=1}^n v_i)v_j=\sum_{i=1}^n \delta_{ij}v_j=v_j$.  Now if we take $e_j\in E^1$ we may use the
equations (2) in the definition of path algebra together with the previous computation to  get $(\sum_{i=1}^n
v_i)e_j=(\sum_{i=1}^n v_i) s(e_j)e_j=s(e_j)e_j=e_j$. In a similar manner we see that  $(\sum_{i=1}^n v_i)e_j^*=e_j^*$.
Since $L(E)$ is generated by $E^0\cup E^1 \cup (E^1)^*$, then it is clear that  $(\sum_{i=1}^n v_i)\alpha=\alpha$ for
every $\alpha\in L(E)$, and analogously $\alpha(\sum_{i=1}^n v_i)=\alpha$ for every  $\alpha\in L(E)$.  Now suppose
that $E^0$ is infinite. Consider a finite subset  $\{a_i\}_{i=1}^t$ of $L(E)$ and use Lemma \monomial\ to write
$a_i=\sum_{s=1}^{n_i}k_s^i v_s^i + \sum_{l=1}^{m_i} c_l^i p_l^i$  where $k_s^i,c_l^i\in K-\{0\}$, and $p_l^i$ are
monomials of type (b). Then with the same ideas as above it is  not difficult to prove that for
$V=\bigcup_{i=1}^t\{v_s^i,s(p_l^i),r(p_l^i):s=1,\dots,n_i;l=1,\dots,m_i\}$, then $\alpha=\sum_{v\in V} v$ is a finite
sum of vertices such that $\alpha a_i=a_i \alpha =a_i$ for every $i$. \endproof

\re{\graded}{\bf Lemma} $L(E)$ is a {\pizarra Z}-graded algebra, with   grading induced by $$deg(v_i)=0 \hbox{ for all
} v_i \in E^0; deg(e_i)=1 \hbox{ and } deg(e_i^*)= -1 \hbox{ for all } e_i \in E^1.$$ That is, $L(E)=\bigoplus_{n\in
\hbox{\pizarra Z}} L(E)_n$, where  $L(E)_0=KE^0 + A_0$, $L(E)_n=A_n$ for $n \neq 0$ where $$A_n=\sum \{ke_{i_1}\dots
e_{i_\s}e_{j_1}^*\dots e_{j_\t}^*:\ \s+\t>0,\  e_{i_s}\in E^1,\ e_{i_t}\in (E^1)^*,\ k\in K,\ \s-\t=n\}.$$ \endre

\proof The fact that $L(E)=\sum_{n\in \hbox{\pizarra Z}} L(E)_n$  follows from Lemma \monomial.  The grading on $L(E)$
follows directly from the fact that $A(\widehat{E})$ is {\pizarra Z}-graded, and that the relations CK1 and CK2 are
homogeneous in this grading.
\endproof

Note that by virtue of Lemma \graded\ we can define the {\bfm degree} of an  arbitrary polynomial in $L(E)$ as the
maximum of the degrees of its monomials. We say that a monomial in $L(E)$ is a {\bfm real path} (resp. a {\bfm ghost
path}) if it contains no terms of the form $e_i^*$ (resp. $e_i$); we say that $p\in L(E)$ is a polynomial in  {\bfm
only real edges} (resp. in  {\bfm only ghost edges}) if it is a sum of real (resp. ghost) paths.

For a path $q=q_1\dots q_n$, we denote by $q^*$ the ghost path $q_n^*\dots q_1^*$. If $\alpha \in L(E)$ and $d \in
{\Bbb Z}^+$, then we say that $\alpha$ is {\bfm representable as an element of degree $d$ in real (resp. ghost) edges}
in case $\alpha$ can be written as a sum of monomials from the spanning set $\{pq^* \mid p,q$ are paths in $E \}$ given
by Lemma \monomial, in such a way that $d$ is the maximum length of a path $p$ (resp. $q$) which appears in such
monomials. We note that an element of $L(E)$ may be representable as an element of different degrees in real (resp.
ghost) edges, depending on the particular representation used for $\alpha$. For instance, for $E$ as in Example
1.4(ii), $xx^{-1}$ is representable as an element of degree $0$ in real edges in $L(E)$, as $xx^{-1} = 1$.

%%%%%%%%%%%%%%%%%%%%%%%%%%%%%%%%%%%%%%%%%%%%%%%%%%%%%%%%%%%%%%%%%%%%%%%%%%%%%%%%
%% %%%%%%%%%%%%%%%  CLOSED PATHS %%%%%%%%%%%%%%%%%%%
%%%%%%%%%%%%%%%%%%%%%%%%%%%%%%%%%%%%%%%%%%%%%%%%%%%%%%%%%%%%%%%%%%%%%%%%%%%%%%%% %%
\bigskip

\parrafo{Closed paths}

\medskip

Certain paths in the graph $E$ will play a central role in the structure of the Leavitt path algebra $L(E)$.

\re{\exit}{}{\bf Definitions.}  An edge $e$ is an {\bfm exit} to the path $\mu=\mu_1 \dots \mu_n$ if there exists $i$
such that $s(e)=s(\mu_i)$ and $e\neq \mu_i$.

A {\bfm closed path based at $v$} is a  path $\mu=\mu_1 \dots \mu_n$, with $\mu_j \in E^1$, $n\geq 1$ and such that
$s(\mu)=r(\mu)=v$. Denote by  $CP(v)$ the set of all such paths.

A {\bfm closed simple path based at $v$}  is a closed path based at $v$, $\mu=\mu_1 \dots \mu_n$, such that
$s(\mu_j)\neq v$ for every $j>1$.  Denote by $CSP(v)$ the set of all such paths.  \endre

{\bf Remark.} Note that a cycle is a closed  simple path based at any of its vertices, but not every closed simple path
based at $v$ is a cycle because a closed simple path may visit some of its vertices (but not $v$) more than
once. Moreover, every closed simple path is in  particular a closed path, while the converse is false.

\re{\productcsp}{\bf Lemma} Let $\mu,\nu\in CSP(v)$.  Then $\mu^*\nu=\delta_{\mu,\nu}v$.
\endre
\proof We first assume $\alpha$ and $\beta$ are arbitrary paths and write $\alpha=e_{i_1}\dots e_{i_\s}$ and
$\beta=e_{j_1}\dots e_{j_\t}$.

Case 1: $deg(\alpha)=deg(\beta)$ but $\alpha\neq \beta$. Define $b\geq 1$ the subindex of the first edge where the
paths $\alpha$ and $\beta$ differ. That is, $e_{i_a}=e_{j_a}$ for every $a<b$ but $e_{i_b}\neq e_{j_b}$. Then
$$\eqalign{ \alpha^*\beta& = e_{i_\s}^*\dots e_{i_1}^* e_{j_1}\dots e_{j_\t}= e_{i_\s}^*\dots e_{i_2}^*
r(e_{j_1}) e_{j_2}\dots e_{j_\t} \cr & = \delta_{r(e_{j_1}),s(e_{j_2})} e_{i_\s}^*\dots e_{i_2}^* e_{j_2}\dots e_{j_\t}
= \dots  \cr & = \delta_{r(e_{j_1}),s(e_{j_2})} \dots \delta_{r(e_{j_{b-1}}),s(e_{j_b})} e_{i_\s}^*\dots e_{i_b}^*
e_{j_b}\dots e_{j_\t}=0. \cr} $$

Case 2: $\alpha=\beta$. Proceeding as above, $\alpha^*\beta=\delta_{r(e_{i_1}),s(e_{i_2})} \dots
\delta_{r(e_{i_{\s-1}}),s(e_{i_\s})} r(e_{i_\s})=r(\alpha)$.

Case 3: Now let $\mu, \nu \in CSP(v)$ with $deg(\mu)<deg(\nu)$. Write $\nu=\nu_1 \nu_2$ where $deg(\nu_1)=deg(\mu),\
deg(\nu_2)>0$. Now if $\mu=\nu_1$ then we have that $v=r(\mu)=r(\nu_1)=s(\nu_2)$, contradicting that $\nu\in CSP(v)$,
so $\mu \neq \nu_1$ and thus case 1 applies to obtain $\mu^*\nu=\mu^*\nu_1\nu_2=0$.

The case $deg(\mu)>deg(\nu)$ is analogous to case 3 by changing the roles of $\mu$ and
$\nu$.
\endproof

\re{\factorcp}{\bf Lemma} For every $p\in CP(v)$ there exist unique $c_1,\dots,c_m\in CSP(v)$ such that $p=c_1\dots
c_m$.
\endre
\proof Write $p=e_{i_1}\dots e_{i_n}$. Let $T=\{t\in \{1,\dots,n\}: r(e_{i_t})=v\}$ and list $t_1<\dots<t_m=n$ all the
elements of $T$. Then $c_1=e_{i_1}\dots e_{i_{t_1}}$ and $c_j=e_{i_{t_{j-1}}}\dots e_{i_{t_j}}$ for $j>1$ give the
desired decomposition.

To prove the uniqueness, write $p=c_1\dots c_r=d_1\dots d_s$ with $c_i,d_j\in CSP(v)$. Multiply by $c_1^*$ on the left
and use Lemma \productcsp\ to obtain $0\neq vc_2\dots c_r=c_1^*d_1\dots d_s$, and therefore by Lemma \productcsp\ again
$c_1=d_1$. Now an induction process finishes the proof.
\endproof

\re{\rd}{}{\bf Definition.} For $p\in CP(v)$ we define the {\bfm return degree (at $v$)} of $p$ to be the number $m\geq
1$ in the decomposition above. (So, in particular, $CSP(v)$ is the subset of $CP(v)$ having return degree equal one.)
We denote it by $RD(p)=RD_{v}(p)=m$. We extend this notion to vertices by setting $RD_{v}(v)=0$, and to nonzero linear
combinations of the form $\sum k_s p_s$, with $p_s \in CP(v)\cup \{v\}$ and $k_s\in K-\{0\}$ by: $RD(\sum k_s p_s)=max
\{RD(p_s)\}$.
\endre

\re{\entries}{\bf Lemma} For a graph $E$ the following conditions are equivalent.
\item{(i)} Every cycle has an exit.
\item{(ii)} Every closed path has an exit.
\item{(iii)} Every closed simple path has an exit.
\item{(iv)} For every $v_i\in E^0$, if $CSP(v_i)\neq \emptyset$, then there exists $c\in CSP(v_i)$ having an exit.
\endre
\proof (ii) $\Rightarrow$ (iii) $\Rightarrow$ (i) is trivial by definition, and (iii) $\Rightarrow$ (iv)
is obvious.

(i) $\Rightarrow$ (ii). Consider $\mu\in CP(v_i)$. First by Lemma \factorcp\ we can factor $\mu=c^{(1)}\dots c^{(m)}$,
where $c^{(j)}\in CSP(v_i)$, and we examine $c^{(m)}$. If it is cycle then we can find an exit for it, and therefore
for $\mu$, by hypothesis. If not, $c^{(m)}$ visits a vertex (different from $v_i$) more than once. Write
$c^{(m)}=c_1^{(m)}\dots c_s^{(m)}$ with each $c_i^{(m)} \in E^1$ and let $c_{s_0}^{(m)}$ be the last edge for which
$s(c_j^{(m)})\in \{s(c_i^{(m)}):1\leq i \leq s, i\neq j \}$. Thus, there exists $s_1<s_0$ such that
$s(c_{s_0}^{(m)})=s(c_{s_1}^{(m)})$. We have several possibilities:

Case 1: $c_{s_0}^{(m)}=c_{s_1}^{(m)}$ and $s_0<s$. Then $r(c_{s_0}^{(m)})=r(c_{s_1}^{(m)})$; that is,
$s(c_{s_0+1}^{(m)})=s(c_{s_1+1}^{(m)})$, which contradicts the choice of $c_{s_0}^{(m)}$.

Case 2: $c_{s_0}^{(m)}=c_{s_1}^{(m)}$ and $s_0=s$. This means that $r(c_{s_1}^{(m)})=r(c_1^{(m)})=v_i$, which is
impossible because $c^{(m)}\in CSP(v_i)$.

Case 3: $c_{s_0}^{(m)}\neq c_{s_1}^{(m)}$. In this case $c_{s_1}^{(m)}$ is an exit for $c^{(m)}$, and then for $\mu$.

In each case we reach a contradiction or we find an exit for $\mu$, as needed.

(iv) $\Rightarrow$ (iii). Consider $c^{(1)}\in CSP(v_i)$. By hypothesis we find $c^{(2)}\in CSP(v_i)$ having an exit.
If $c^{(1)}=c^{(2)}$ we are done. If not, we write $c^{(1)}=e_{i_1}\dots e_{i_s}$, $c^{(2)}=e_{j_1}\dots e_{j_r}$ and
proceed by steps:

Step 1: If $e_{i_1}\neq e_{j_1}$, since $s(e_{i_1})=s(e_{j_1})=v_i$, then $e_{j_1}$ is an exit for $c^{(1)}$.

Step 2: If $e_{i_1}=e_{j_1}$ then $r(e_{i_1})=r(e_{j_1})$; that is, $s(e_{i_2})=s(e_{j_2})$.

Step 3: If $e_{i_2}\neq e_{j_2}$, then as in Step 1, $e_{j_2}$ is an exit for $c^{(1)}$.

Step 4: If $e_{i_2}=e_{j_2}$, then continue as in Step 2.

With this process, we either find an exit or we run out of edges in one path but not in the other (because $c^{(1)}\neq
c^{(2)}$). Thus:

Case 1: $c^{(1)}=c^{(2)}e_{i_t}\dots e_{i_s}$ for $t\leq s$. But this is impossible because $s(e_{i_t})=r(c^{(2)})=v_i$
and $c^{(1)}\in CSP(v_i)$.

Case 2: $c^{(2)}=c^{(1)}e_{j_q}\dots e_{j_r}$ for $q\leq r$, which is similarly impossible.

In any case, we reach a contradiction or we are able to find an exit for $c^{(1)}$, and this finishes the proof.
\endproof

%%%%%%%%%%%%%%%%%%%%%%%%%%%%%%%%%%%%%%%%%%%%%%%%%%%%%%%%%%%%%%%%%%%%%%%%%%%%%%%%
%%%%%%%%%%%%%%%  SIMPLICITY OF L(E)
%%%%%%%%%%%%%%%%%%%
%%%%%%%%%%%%%%%%%%%%%%%%%%%%%%%%%%%%%%%%%%%%%%%%%%%%%%%%%%%%%%%%%%%%%%%%%%%%%%%%
\bigskip

\parrafo{Simplicity of $L(E)$}

\medskip

In this final section we build the algebraic machinery necessary to obtain our main result, Theorem (3.11).

\re{\decreasedegreereal}{\bf Proposition} Let $E$ be a graph with the property that every cycle has an exit. If $\alpha
\in L(E)$ is a polynomial in only real edges with $deg(\alpha)>0$, then there exist $a,b\in L(E)$ such that $a\alpha
b\neq 0$ is a polynomial in only real edges and $deg(a\alpha b)<deg(\alpha)$.
\endre
\proof Write $\alpha=\sum_{e_i \in E^1} e_i \alpha_{e_i}+\sum_{v_l \in E^0} k_l v_l$, where $\alpha_{e_i}$ are
polynomials in only real edges,  and $deg(\alpha_{e_i})<deg(\alpha)=m$.

Case (A): $k_l=0$ for every $l$. Since $\alpha \neq 0$, there exists $i_0$ such that $e_{i_0} \alpha_{e_{i_0}}\neq 0$.
Let $b\in L(E)$ have $\alpha b = \alpha$; such exists by Lemma \unital. Then $a=e_{i_0}^*$, $b$ give $e_{i_0}^* \alpha
b=\alpha_{e_{i_0}}\neq 0$ is a polynomial in only real edges and $deg(\alpha_{e_{i_0}})<deg(\alpha)$.

Case (B): There exists $k_{l_0}\neq 0$. Then we can write $$v_{l_0}\alpha v_{l_0}=k_{l_0}v_{l_0}+\sum_{p\in
CP(v_{l_0})} k_p p,\ k_p\in K.$$ Note that this is a polynomial in only real edges, and is nonzero because
$k_{l_0}$ is nonzero.

Case (B.1): $deg(v_{l_0}\alpha v_{l_0})<deg(\alpha)$. Then we are done with $a=v_{l_0}$ and $b=v_{l_0}$.

Case (B.2): $deg(v_{l_0}\alpha v_{l_0})=deg(\alpha)=m>0$. Then there exists $p_0\in CP(v_{l_0})$ such that $k_{p_0}
p_0\neq 0$. Now by Lemma \factorcp, we can write $p_0=c_1\dots c_\s$, $\s\geq 1$ and thus $CSP(v_{l_0})\neq \emptyset$.
We apply now Lemma \entries\ to find $c_{s_0}\in CSP(v_{l_0})$ which has $e_{i_0}$ as an exit, that is, if
$c_{s_0}=e_{i_1}\dots e_{i_{s_0}}$ then there exists $j\in \{1,\dots,s_0\}$ such that $s(e_{i_j})=s(e_{i_0})$ but
$e_{i_j}\neq e_{i_0}$. Since $s(e_{i_j})=s(e_{i_0})$ we can therefore build the path given by $z=e_{i_1}\dots
e_{i_{j-1}}e_{i_0}$. This path has $c_{s_0}^*z=0$ because $c_{s_0}^*z=e_{i_{s_0}}^*\dots e_{i_1}^*e_{i_1}\dots
e_{i_{j-1}}e_{i_0}= \dots = e_{i_{s_0}}^*\dots e_{i_j}^*e_{i_0}=0$. (We will use this observation later on.) Again
Lemma \factorcp\ allows us to write $$v_{l_0}\alpha v_{l_0}=k_{l_0}v_{l_0}+ \sum_{c_s\in CSP(v_{l_0})} c_s
\alpha_{c_s}^{(1)}, \tag{\dag}$$ where $\gamma=RD(v_{l_0}\alpha v_{l_0})>0$, and $\alpha_{c_s}^{(1)}$ are polynomials
in only real edges satisfying $RD(\alpha_{c_s}^{(1)})<\gamma$.

We now present a process in which we decrease the return degree of the polynomials by multiplying on both sides by
appropriate elements in $L(E)$. In the sequel we will often make use of Lemma \productcsp\ without mentioning it
explicitly. In particular, multiplying $(\dag)$ on the left by $c^*_{s_0}$ gives $$c_{s_0}^*(v_{l_0}\alpha
v_{l_0})=k_{l_0}c_{s_0}^*+ \alpha_{c_{s_0}}^{(1)}. \tag{\ddag}$$

Case 1: $\alpha_{c_{s_0}}^{(1)}=0$. Then $A=c_{s_0}^*$ and $B=c_{s_0}$ are such that $A(v_{l_0}\alpha
v_{l_0})B=k_{l_0}v_{l_0}\neq 0$ is a polynomial in only real edges and $RD(A(v_{l_0}\alpha
v_{l_0})B)=0<\gamma=RD(v_{l_0}\alpha v_{l_0})$.

Case 2: $\alpha_{c_{s_0}}^{(1)}\neq 0$ but $RD(\alpha_{c_{s_0}}^{(1)})=0$. Then $\alpha_{c_{s_0}}^{(1)}=k^{(2)}v_{l_0}$
for some $0 \neq k^{(2)} \in K$. Using the path $z$ with an exit for $c_{s_0}^*$ we have: $z^*c_{s_0}^*(v_{l_0}\alpha
v_{l_0})z=z^*(k_{l_0}c_{s_0}^*+k^{(2)}v_{l_0})z=z^*(0+k^{(2)}z)=k^{(2)}r(z)\neq 0$. So we have $A=z^*c_{s_0}^*$ and
$B=z$ such that $A(v_{l_0}\alpha v_{l_0})B\neq 0$ is a polynomial in only real edges and $RD(A(v_{l_0}\alpha
v_{l_0})B)=0<\gamma=RD(v_{l_0}\alpha v_{l_0})$.

Case 3: $RD(\alpha_{c_{s_0}}^{(1)})>0$. We can write $$\alpha_{c_{s_0}}^{(1)}=k^{(2)}v_{l_0}+ \sum_{c_s\in
CSP(v_{l_0})} c_s \alpha_{c_s}^{(2)},$$ where $\alpha_{c_s}^{(2)}$ are polynomials in only real edges with return
degree less than the return degree of $\alpha_{c_{s_0}}^{(1)}$. Now $0<RD(\alpha_{c_{s_0}}^{(1)})<\gamma$ implies
$\gamma\geq 2$. Multiply $({\ddag})$ by $c_{s_0}^*$ to get $$(c_{s_0}^*)^2(v_{l_0}\alpha v_{l_0})=k_{l_0}(c_{s_0}^*)^2+
k^{(2)}c_{s_0}^* + \alpha_{c_{s_0}}^{(2)}. \tag{\S}$$

We are now in position to proceed in a manner analogous to that described in Cases 1, 2, and 3 above.

\medskip

Case 3.1: $\alpha_{c_{s_0}}^{(2)}=0$. Then $(c_{s_0}^*)^2(v_{l_0}\alpha
v_{l_0})(c_{s_0})^2=k_{l_0}v_{l_0}+k^{(2)}c_{s_0}$ and hence we have found $A=(c_{s_0}^*)^2$ and $B=(c_{s_0})^2$ such
that $A(v_{l_0}\alpha v_{l_0})B\neq 0$ is a polynomial in only real edges and $RD(A(v_{l_0}\alpha v_{l_0}) B)=1<2\leq
\gamma =RD(v_{l_0}\alpha v_{l_0})$.

Case 3.2: $\alpha_{c_{s_0}}^{(2)}\neq 0$ but $RD(\alpha_{c_{s_0}}^{(2)})=0$. Then
$\alpha_{c_{s_0}}^{(2)}=k^{(3)}v_{l_0}$ for some $0 \neq k^{(3)} \in K$, and then $z^*(c_{s_0}^*)^2(v_{l_0}\alpha
v_{l_0})z=z^*(k_{l_0}(c_{s_0}^*)^2+k^{(2)}c_{s_0}^*+k^{(3)}v_{l_0})z=z^*(0+k^{(3)}z)=k^{(3)}r(z)\neq 0$. Thus, we get
$A=z^*(c_{s_0}^*)^2$ and $B=z$ such that $A(v_{l_0}\alpha v_{l_0})B\neq 0$ is a polynomial in only real edges and
$RD(A(v_{l_0}\alpha v_{l_0})B)=0<\gamma=RD(v_{l_0}\alpha v_{l_0})$.

Case 3.3: $RD(\alpha_{c_{s_0}}^{(2)})>0$. We write $$\alpha_{c_{s_0}}^{(2)}=k^{(3)}v_{l_0}+ \sum_{c_s\in CSP(v_{l_0})}
c_s \alpha_{c_s}^{(3)},$$ where $\alpha_{c_s}^{(3)}$ are polynomials in only real edges with return degree less than
the return degree of $\alpha_{c_{s_0}}^{(2)}$. Now $0<RD(\alpha_{c_{s_0}}^{(2)})<RD(\alpha_{c_{s_0}}^{(1)})<\gamma$
implies $\gamma\geq 3$. And by multiplying $({\S})$ by $c_{s_0}^*$ we get $(c_{s_0}^*)^3(v_{l_0}\alpha
v_{l_0})=k_{l_0}(c_{s_0}^*)^3+ k^{(2)}(c_{s_0}^*)^2+ k^{(3)}c_{s_0}^* + \alpha_{c_{s_0}}^{(3)}$.

We continue the process of analyzing each such equation by considering three cases.  If at any stage either of the
first two cases arise, we are done.  But since at each stage the third case can occur only by producing elements of
subsequently smaller return degree, then after at most $\gamma$ stages we must have one of the first two cases.

Thus, by repeating this process at most $\gamma$ times we are guaranteed to find $\widetilde{A},\widetilde{B}$ such
that $\widetilde{A}(v_{l_0}\alpha v_{l_0})\widetilde{B}\neq 0$ is a polynomial in only real edges and
$RD(\widetilde{A}(v_{l_0}\alpha v_{l_0})\widetilde{B})=0$. But this then gives $0=deg(\widetilde{A}(v_{l_0}\alpha
v_{l_0})\widetilde{B})<deg(\alpha)$. So $a=\widetilde{A}v_{l_0}$ and $b=v_{l_0}\widetilde{B}$ are the desired elements.
\endproof

\re{\degreezero}{\bf Corollary} Let $E$ be a graph with the property that every cycle has an exit. If $\alpha\neq 0$ is
a polynomial in only real edges then there exist $a,b\in L(E)$ such that $a\alpha b\in E^0$.
\endre
\proof Apply Proposition \decreasedegreereal\ as many times as needed ($deg(\alpha)$ at most) to find $a',b'$ such that
$a'\alpha b'$ is a nonzero polynomial in only real edges with $deg(a'\alpha b')=0$; that is, $a'\alpha b'=\sum_{i=1}^t
k_i v_i \neq 0$.  So there exists $j$ with $k_{j}\neq 0$, and finally $a=k_{j}^{-1}a'$ and $b=b'v_{j}$ give that
$a\alpha b=v_{j}\in E^0$.
\endproof

\re{\idealJ}{\bf Corollary} Let $E$ be a graph with the property that every cycle has an exit. If $J$ is a ideal of
$L(E)$ and contains a nonzero polynomial in only real edges, then $E^0\cap J\neq \emptyset$.
\endre
\proof Straightforward by Corollary \degreezero.
\endproof

In order to extend all the previous results of this section to analogous results about polynomials in only ghost edges,
 we define an involution
in $L(E)$.

\re{\involution}{\bf Lemma} $L(E)$ can be equipped with an involution $x\mapsto\overline{x}$ defined in the monomials
by:
\item{(a)} $\overline{k_iv_i}=k_iv_i$ with $k_i\in K$ and $v_i\in E^0$,
\item{(b)} $\overline{ke_{i_1}\dots e_{i_\s}e_{j_1}^*\dots e_{j_\t}^*}=ke_{j_\t}\dots e_{j_1}e_{i_\s}^*\dots e_{i_1}^*$
where $k\in K;\ \s,\ \t\geq 0,\ \s+\t>0,\ e_{i_s}\in E^1$ and $e_{j_t}\in (E^1)^*$,

and extending linearly to $L(E)$.
\endre
\proof The proposed map is well defined by Lemma \monomial, and it is linear by definition. It is easily shown to
satisfy $\overline{xy}=\overline{y}\ \overline{x}$ and $\overline{\overline{x}}=x$ for every $x,y\in L(E)$. It is also
straightforward to check that the map is compatible with the relations defining $L(E)$.
\endproof

\re{\generalizeghost}{}{\bf Remark.} Note that the involution transforms a polynomial in only real edges into a
polynomial in only ghost edges and vice versa. If $J$ is an ideal of $L(E)$ then so is $\overline{J}$. We note here
that while Leavitt path algebras behave somewhat like their $C^*$-algebra siblings, they are indeed different in many
respects.    For instance, whereas in $C^*$-algebras every two-sided ideal $J$ is self-adjoint (i.e. $\overline{J} =
J$), this is not the case in the Leavitt path algebras setting. For instance, let $L(E) = K[x,x^{-1}]$ as in Example
\examples\ (ii), and let $J$ be the ideal $<1+x+x^3>$ of $L(E)$. Then $J$ is not self-adjoint, as follows: if
$\overline{J}=J$, then $f(x) = 1 + x^{-1} + x^{-3} \in J$ and thus $x^3f(x) = 1 + x^2 + x^3 \in J$. Now $K[x,x^{-1}]$
being a unital commutative ring implies that there exists $p=\sum_{i=-\infty}^{\infty} a_ix^i$ with
$p(1+x+x^3)=1+x^2+x^3$. A degree argument on the highest power on the left hand side of the previous equation leads to
$a_i=0$ for every $i\geq 1$. By reasoning in a similar fashion on the lowest power we also get $a_i=0$ for every $i\leq
-1$, that is, $p=a_0$, which is absurd.
\endre

We can define sets and quantities for ghost paths analogous to those given for real paths. Using the involution given
in Lemma \involution\ we can then analogously prove the following three results.

\re{\increasedegreeghost}{\bf Proposition} Let $E$ be a graph with the property that every cycle has an exit. If
$\alpha \in L(E)$ is a polynomial in only ghost edges with $deg(\overline{\alpha})>0$ then there exist $a,b\in L(E)$
such that $a\alpha b\neq 0$ is a polynomial in only ghost edges and $deg(\overline{a\alpha b})<deg(\overline{\alpha})$.
\endre

\re{\degreezeroghost}{\bf Corollary} Let $E$ be a graph with the property that every cycle has an exit. If $\alpha\neq
0$ is a polynomial in only ghost edges then there exist $a,b\in L(E)$ such that $a\alpha b\in E^0$.
\endre

\re{\idealJghost}{\bf Corollary} Let $E$ be a graph with the property that every cycle has an exit. If $J$ is an ideal
of $L(E)$ and contains a nonzero polynomial in only ghost edges, then $E^0\cap J\neq \emptyset$.
\endre

For a graph $E$ we define a preorder $\leq$ on the vertex set $E^0$ given by:
$$v\leq w \hbox{ if and only if }v=w\hbox{ or there is a path } \mu \hbox{ such that } s(\mu)=v \hbox{ and } r(\mu)=w.$$
We say that a subset $H\subseteq E^0$ is {\bfm hereditary} if $w\in H$ and $w\leq v$ imply $v\in H$. We say that $H$ is
{\bfm saturated} if whenever $s^{-1}(v)\neq \emptyset$ and $\{r(e):s(e)=v\}\subseteq H$, then $v\in H$.  (In other
words, $H$ is saturated if, for any vertex $v$ in $E$, if {\it all} of the range vertices $r(e)$ for those edges $e$
having $s(e)=v$ are in $H$, then $v$ must be in $H$ as well.)

\re{\hereditarysaturated}{\bf Lemma} If $J$ is an ideal of $L(E)$, then $J\cap E^0$ is a hereditary and saturated
subset of $E^0$.
\endre
\proof We first show that $J\cap E^0$ is hereditary.  Consider $v,w \in E^0$ such that $v\in J$ and $v\leq w$.
 By the
definition of the preorder we can find a path $\mu=\mu_1\dots \mu_n$ such that $s(\mu_1)=v$ and $r(\mu_n)=w$. Apply
that $J$ is an ideal to get that $\mu_1^*v\mu_1=\mu_1^*\mu_1=r(\mu_1)=s(\mu_2)\in J$. Repeating this argument $n$
times, we get that $r(\mu_n)=w\in J$.

Now we see that $J\cap E^0$ is saturated: consider a vertex $v$ with $s^{-1}(v)\neq \emptyset$ and
$\{r(e):s(e)=v\}\subseteq J$. The first condition implies that $v$ is not a sink, so CK2 applies and we obtain
$v=\sum_{\{e_j\in E^1:s(e_j)=v\}}e_je_j^*$. If we take $e_j$ such that $s(e_j)=v$, then by hypothesis we have that
$r(e_j)\in J$ and therefore $e_j=e_jr(e_j)\in J$. Now applying CK2 we conclude that $v\in J$.
\endproof

\re{\corollary}{\bf Corollary} Let $E$ be a graph with the following properties:
\item{(i)} The only hereditary and saturated subsets of $E^0$ are $\emptyset$ and $E^0$.
\item{(ii)} Every cycle has an exit.

If $J$ is a nonzero ideal of $L(E)$ which contains a polynomial in only real edges (or a polynomial in only ghost
edges), then $J=L(E)$.
\endre
\proof Apply Corollaries \idealJ\ or \idealJghost\ to get that $J\cap E^0\neq \emptyset$. Now by Lemma
\hereditarysaturated\ and (i) we have  $J\cap E^0=E^0$. Therefore $J$ contains a set of local units by Lemma \unital,
and hence $J=L(E)$.
\endproof

We are now in position to prove the main result of this article.

\re{\simple}{\bf Theorem} Let $E$ be a row-finite graph. Then the Leavitt path algebra $L(E)$ is simple if and only if
$E$ satisfies the following conditions.
\item{(i)} The only hereditary and saturated subsets of $E^0$ are $\emptyset$
and $E^0$, and
\item{(ii)} Every cycle in $E$ has an exit.
\endre
\proof First we assume that (i) and (ii) hold and we will show that $L(E)$ is simple. Suppose that $J$ is a nonzero
ideal of $L(E)$.  Choose $0\neq \alpha \in J$ representable as an element having minimal degree in the real edges. If
this minimal degree is $0$, then $\alpha$ is a polynomial in only ghost edges, so that by Corollary \corollary\ we have
$J=L(E)$. So suppose this degree in real edges is at least 1.  Then we can write $$\alpha =\sum_{n=1}^m e_{i_n}
\alpha_{e_{i_n}} + \beta$$ where $m\geq 1$, $e_{i_n}\alpha_{e_{i_n}}\neq 0$ for every $n$, and each $\alpha_{e_{i_n}}$
is representable as an element of degree less than that of $\alpha$ is real edges, and $\beta$ is a polynomial in only
ghost edges (possibly zero).

Suppose $v$ is a sink in $E$. Then we may assume $v\beta =0$, as follows. Multiplying the displayed equation by $v$ on
the left gives $v\alpha = v\sum_{n=1}^m e_{i_n} \alpha_{e_{i_n}} + v\beta$.  But since $v$ is a sink we have $v e_{i_n}
=0$ for all $1 \leq n \leq m$, so that $v\alpha = v\beta \in J$.  But $v\beta \neq 0$ would then yield a nonzero
element of $J$ in only ghost edges, so that again by Corollary \corollary\ we have $J=L(E)$.

For an arbitrary edge $e_j\in E^1$, we have two cases:

Case 1: $j\in \{i_1,\dots, i_m\}$. Then $e_j^*\alpha=\alpha_{e_j}+e_j^*\beta\in J$. If this element is nonzero it would
be representable as an element with smaller degree in the real edges than that of $\alpha$, contrary to our choice. So
it must be zero, and hence $\alpha_{e_j}=-e_j^*\beta$, so that $e_j\alpha_{e_j}=-e_je_j^*\beta$.

Case 2: $j\not\in \{i_1,\dots, i_m\}$. Then $e_j^*\alpha=e_j^*\beta\in J$.  If $e_j^*\beta \neq 0$,
then as before we would have a nonzero element of $J$ in only ghost edges, so that $J=L(E)$ and we are done.
So we may assume that $e_j^*\beta = 0$, so that in particular we have $0=-e_je_j^*\beta$.

Now let $S_1=\{v_j \in E^0:v_j = s(e_{i_n}) \hbox{ for some } 1\leq n \leq m\}$, and let $S_2 =
\{v_{k_1},...,v_{k_t}\}$ where $(\sum_{i=1}^t v_{k_i})\beta = \beta$.  (Such a set $S_2$ exists by Lemma \unital.) We
note that $w\beta = 0$ for every $w\in E^0 - S_2$.  Also, by definition there are no sinks in $S_1$, and by a previous
observation we may assume that there are no sinks in $S_2$.     Let $S = S_1 \cup S_2$. Then in particular we have
$(\sum _{v\in S} v)\beta = \beta$.

We now argue that in this situation $\alpha$ must be zero, which will contradict our original choice of $\alpha$ and
thereby complete the proof.  To this end,
 $$\eqalign{
 \alpha & =\sum_{n=1}^m e_{i_n} \alpha_{e_{i_n}} + \beta = \sum_{n=1}^m -e_{i_n}e_{i_n}^*\beta + \beta \hskip5mm \hbox{
  (by Case 1)} \cr
 & = \sum_{n=1}^m -e_{i_n}e_{i_n}^*\beta - (\sum_{j\notin \{i_1,...,i_m\}, s(e_j)\in S} e_je_j^*) \beta  + \beta  \cr
 & \hskip2cm \hbox{ (by Case 2, the newly subtracted terms equal 0) }\cr
 & =-(\sum _{v\in S} v)\beta + \beta  \hskip5mm \hbox{ (no sinks in } S \hbox{ implies that CK2 applies at each }v\in S)
\cr
 & = -\beta + \beta = 0.\cr}
 $$

Thus we have shown that if $E$ satisfies the two indicated properties, then $L(E)$ is simple.   For the converse, first
suppose that there is a cycle $p$ having no exit. We will prove that $L(E)$ cannot be simple. Let $v$ be the base of
that cycle. We will show that for $\alpha=v+p$, $<\alpha>$ is a nontrivial ideal of $L(E)$ because $v\not\in <\alpha>$.
Write $p=e_{i_1}\dots e_{i_\s}$. Since this cycle does not have an exit, for every $e_{i_j}$ there is no edge with
source $s(e_{i_j})$ other than $e_{i_j}$ itself, so that the CK2 relation at this vertex yields
$s(e_{i_j})=e_{i_j}e_{i_j}^*$. This easily implies $pp^*=v$ (we recall here that $p^*p=v$ always holds), and that
$CSP(v)=\{p\}$.

Now suppose that $v\in <\alpha>$.  So there exist nonzero monic monomials $a_n,b_n \in L(E)$ and $c_n\in K$ with
$v=\sum_{n=1}^m c_n a_n\alpha b_n$. Since $v\alpha v = \alpha$, by multiplying by $v$ if necessary we may assume that
$va_n v=a_n$ and $vb_n v=b_n$ for all $1\leq n \leq m$.

We claim that for each $a_n$ (resp. $b_n$) there exists an integer $u(a_n)\geq 0$ (resp. $u(b_n)\geq 0$) such that $a_n
= p^{u(a_n)}$ or  $a_n = (p^*)^{u(a_n)}$ (resp. $b_n = p^{u(b_n)}$ or  $b_n = (p^*)^{u(b_n)}$).

Now $a_1$ is of the form $e_{k_1}\dots e_{k_c}e_{j_1}^*\dots e_{j_d}^*$ with $c,d\geq 1$. (Otherwise we are in a simple
case that will be contained in what follows.) Since $a_1$ starts and ends in $v$ we can consider the elements:
$g=min\{z:r(e_{j_z}^*)=v\}$ and $f=max\{z:s(e_{k_z})=v\}$, and we will focus on $a'_1=e_{k_f}\dots
e_{k_c}e_{j_1}^*\dots e_{j_g}^*$.

First, since $v=r(e_{j_g}^*)=s(e_{j_g})$ and $e_{i_1}$ is the only edge coming from $v$, then $e_{j_g}=e_{i_1}$. Now,
$s(e_{j_{g-1}})=r(e_{j_{g-1}}^*)=s(e_{j_g}^*)=r(e_{j_g})=r(e_{i_1})=s(e_{i_2})$, and again the only edge coming from
$s(e_{i_2})$ is $e_{i_2}$ and therefore $e_{j_{g-1}}=e_{i_2}$. This process must stop before we run out of edges of $p$
because by our choice of $g$ we have that $v\not\in \{r(e_{j_z}^*):z<g\}$. So in the end there exists $\gamma<\s$ such
that $e_{j_1}^*\dots e_{j_g}^*=e_{i_\gamma}^*\dots e_{i_1}^*$.

With the same (reversed) ideas in the paragraph above we can find $\delta<\s$ such that $e_{k_f}\dots
e_{k_c}=e_{i_1}\dots e_{i_\delta}$. Thus, $a'_1=e_{i_1}\dots e_{i_\delta}e_{i_\gamma}^*\dots e_{i_1}^*$, and we
have two cases:

Case 1: $\delta\neq \gamma$. We know that $p$ is a cycle, so that $r(e_{i_\delta})\neq
r(e_{i_\gamma})=s(e_{i_\gamma}^*)$, so $e_{i_\delta}e_{i_\gamma}^*=0$, which is absurd because $a_1\neq 0$.

Case 2: $\delta= \gamma$. In this case $a'_1=p_0p_0^*$ for a certain subpath $p_0$ of $p$, and by using again the
argument of the CK2 relation in this case, we obtain $p_0p_0^*=v$.

Hence, we get $a_1=e_{k_1}\dots e_{k_{f-1}}e_{j_{g+1}}^*\dots e_{j_d}^*=xy^*$, with $x,y\in CP(v)$. (Obviously, the
case $c\geq 1, d=0$ yields $a_1=x$, the case $c=0, d\geq 1$ yields $a_1=y^*$ and $c=d=0$ yields $a_1=v$.) Using Lemma
\factorcp\ we have $x=c^{(1)}\dots c^{(\nu)}$ for some $c^{(\mu)}\in CSP(v)=\{p\}$, and the same happens with $y$. In
this way we have $a_1=p^u(p^*)^v$ for some $u,v\geq 0$, and taking into account that $pp^*=v$ we finally obtain that
$a_1$ is of the form $p^u$ or $(p^*)^u$ for some $u\geq 0$ as claimed.   An identical argument holds for the other
coefficients $a_n$ and $b_n$.

Now since both $p$ and $p^*$ commute with $p,p^*$ and $\alpha$, we use the conclusion of the previous paragraph to
write the sum $v=\sum_{n=1}^m c_n a_n\alpha b_n$ as $v=\alpha P(p,p^*)$ for some polynomial $P$ having coefficients in
$K$. Specifically, $P(p,p^*)$ can be written as $P(p,p^*)=k_{-m}(p^*)^m+\dots+k_0 v+\dots+k_n p^n\in \bigoplus_{j=-m}^n
L(E)_{\s j}$, where $m,n\geq 0$. First, we claim that $k_{-i}=0$ for every $i>0$, as follows.  If not, let $m_0$ be the
maximum $i$ having $k_{-i}\neq 0$. Then $\alpha P(p,p^*)=k_{-m_0}(p^*)^{m_0}+\text{ terms of greater degree }=v$, and
since $m_0>0$ we get that $k_{-m_0}=0$, which is absurd. In a similar way we obtain $k_i=0$ for every $i>0$, and
therefore $P(p,p^*)=k_0 v$. But this would yield $v =\alpha P(p,p^*)=\alpha k_0 v = k_0 \alpha$, which is impossible.

Thus we have shown that if $E$ contains a cycle which has no exit, then $L(E)$ is not simple. Now we will consider the
situation where $E^0$ contains a nontrivial hereditary and saturated subset $H$, and conclude in this case as well that
$L(E)$ is not simple.  To do so, we construct a new graph
$F=(F^0,F^1,r_F,s_F)=(E^0-H,r^{-1}(E^0-H),r|_{E^0-H},s|_{E^0-H})$.  In other words, $F$ is the graph consisting of all
vertices not in $H$, together with all edges whose range is not in $H$. To ensure that $F$ is well-defined, we must
check that $s_F(F^1)\cup r_F(F^1) \subseteq F^0$. That $r_F(F^1) \subseteq F^0$ is evident. On the other hand, if $e\in
F^1$ then $s(e)\in F^0$, since otherwise we have $s(e)\in H$; but since $r(e)\geq s(e)$ and $H$ is hereditary, we get
$r(e)\in H$, which contradicts $e\in F^1$.  So $F$ is a well defined graph.

We now produce a $K$-algebra homomorphism $\Psi:L(E)\to L(F)$.  To do so, we define $\Phi$ on the generators of the
free $K$-algebra $B = K[E^0\cup E^1\cup (E^1)^*]$ by setting $\Phi(v_i)={\chi}_{F^0}(v_i) v_i$,
$\Phi(e_i)=\chi_{F^1}(e_i)e_i$ and $\Phi(e^*_i)=\chi_{(F^1)^*}(e^*_i)e^*_i$ (where $\chi_X$ denotes the usual
characteristic function of a set $X$), and extending to $B$.  In order to factor $\Phi$ through $A(\widehat{E})$ we
need to check that $$<\{v_iv_j-\delta_{ij}v_i:v_i,v_j\in E^0\}\cup\{e_i-e_ir(e_i),e_i-s(e_i)e_i: e_i\in
\widehat{E}^1\}>\hbox{ }\subseteq Ker(\Phi).$$ This is a straightforward computation done by cases, with the only
nontrivial situation arising when $e_i\in F^1$. But then $r(e_i)\not\in H$, and therefore
$\Phi(e_i-e_ir(e_i))=e_i-e_ir(e_i)=0$ in $L(F)$. Now, since $s(e_i)\leq r(e_i)\not\in H$ and $H$ is hereditary then
$s(e_i)\not\in H$, so that $\Phi(e_i-s(e_i)e_i)=e_i-s(e_i)e_i=0$ in $L(F)$.

Now to produce the desired ring homomorphism $\Psi: L(E)\to L(F)$ we need only check that
$\Phi$ factors through the relations ideal
$$<\{e_i^*e_j-\delta_{ij}r(e_j): e_j\in E^1,e_i^*\in (E^1)^*\}\cup\{v_i-\sum_{\{e_j\in E^1:s(e_j)=v_i\}}e_je_j^*:
v_i\in s(E^0)\}>$$ of $A(\widehat{E})$.  That $\Phi(e_i^*e_j-\delta_{ij}r(e_j))=0$ in $L(F)$ is straightforward.  So
now consider $v_i\in s(E^0)$; i.e., consider a vertex $v_i$ which is not a sink in $E$.

Case 1: Suppose $v_i\in H$. Then for every $e_j\in E^1$ with $s(e_j)=v_i$ we have that $e_i\not\in F^1$ (otherwise
$e_i\in F^1$ implies $r(e_i)\not\in H$ and by hereditariness $s(e_j)=v_i\not\in H$). So, $\Phi(v_i-\sum_{\{e_j\in
E^1:s(e_j)=v_i\}}e_je_j^*)=0-\sum_{\{e_j\in E^1:s(e_j)=v_i\}}0\cdot 0=0$.

Case 2: Suppose $v_i\not\in H$ and $v_i\not\in s(F^1)$. Since $v_i\in s(E^0)$ we have $s^{-1}(v_i)\neq \emptyset$. But
since $H$ is saturated there must exist $e_i\in E^1$ such that $s(e_i)=v_i$, but $r(e_i)\not\in H$. That means $e_i\in
F^1$ with $s(e_i)=v_i$, which contradicts the hypothesis that $v_i\not\in s(F^1)$. Thus the saturated condition on $H$
implies that Case 2 configuration cannot occur.

Case 3: Suppose $v_i\not\in H$ but $v_i\in s(F^1)$. Then we have a CK2 relation in $L(F)$ at $v_i$: $$
v_i=\sum_{\{e_j\in F^1:s(e_j)=v_i\}}e_je_j^*.$$ Consider $e_j\in E^1$ such that $s(e_j)=v_i$. If $e_j\in F^1$ then
$\Phi(e_je^*_j)=e_je^*_j$. If $e_j\not\in F^1$ then $\Phi(e_je^*_j)=0$. Thus we get $\Phi(v_i-\sum_{\{e_j\in
E^1:s(e_j)=v_i\}}e_je_j^*)=v_i-\sum_{\{e_j\in F^1:s(e_j)=v_i\}}e_je_j^*$ $ = 0$ by the displayed equation.

Thus we have shown that there exists a $K$-algebra homomorphism $\Psi:L(E)\to L(F)$.
Now consider $Ker(\Psi)\trianglelefteq L(E)$. Since $H\neq \emptyset$ there exists $v\in H$, so $0\neq v\in Ker(\Psi)$.
 Since $H\neq E^0$ there exists $w\in E^0-H$ and in this case
$\Psi(w)=w\neq 0$ so $\Psi \neq 0$. In other words, $0\neq Ker(\Psi) \neq L(E)$, so that $L(E)$ is not
simple.

Thus we conclude that the negation of either condition (i) or condition (ii) yields that $L(E)$ is not simple, which
completes the proof of the theorem.  \endproof

\re{\algorithm}{}{\bf Remark.} If we start with a finite and row-finite graph $E=(E^0,E^1,r,s)$ with
$E^0=\{v_1,\dots,v_n\},E^1=\{e_1,\dots,e_m\}$, there exist algorithms that decide, in a finite number of steps, whether
or not the graph satisfies conditions (i) and/or (ii), and therefore whether or not $L(E)$ is simple. \endre

\re{\corollary}{}{\bf Corollary.} We re-establish the simplicity (or non-simplicity) of the algebras given in Examples
\examples\ above.

\item{(i)} Matrix algebras $M_n(K)$:  Since there are clearly no cycles in $E$,
we need only verify condition (i) in Theorem \simple.  To this end, let $H\neq \emptyset$ be a set of vertices which is
hereditary and saturated. Pick $v_i\in H$. By hereditariness we have that $v_{i+1},\dots,v_n\in H$. Now if we use the
condition of being saturated at $v_{i-1}$ we get that $v_{i-1}\in H$, and inductively $v_{i-1},\dots,v_1\in H$ and
therefore $H=E^0$. Hence there are no nontrivial hereditary and saturated subsets of $E^0$, and Theorem \simple\
applies to give that $M_n(K)=L(E)$ is simple.

\item{(ii)} Laurent polynomial algebras $K[x,x^{-1}]$:  The cycle $x$ does not have an exit,
so by Theorem \simple\ $L(E) \cong K[x,x^{-1}]$ is not simple.  (Indeed, similar to the argument which arises in the
proof of Theorem \simple, it is easy to show that $1 \notin <1+x>$.)

\item{(iii)} Leavitt algebras $L(1,n)$ for $n\geq 2$: The
conditions in Theorem \simple\ are clearly satisfied here, so $L(1,n)$ is simple, as was established in \cite{\Ltwo,
Theorem 2}.

\endre

\re{\CnExample}{}{\bf Example.}  Let $C_n$ denote the graph having $n$ vertices and $n$ edges, where the edges form a
single cycle.  (In particular, the graph described in Example \examples\ (ii) is the graph $C_1$.)  Then $L(C_n)$ is
not simple for all $n$, since the single cycle contains no exit.

\endre

\re{\CKExample}{}{\bf Example.}  The Cuntz-Krieger algebra ${\Cal C}{\Cal K}_A(K)$ of a finite matrix $A$ is defined in
\cite{\AGGP, example 2.5}.  For a finite graph $E$ we can define the {\bfm edge matrix} $A_E$ associated to $E$; $A_E$
is the $n \times n$ matrix with entries $a_{ij}=\delta_{r(e_i),s(e_j)}$, where $n=|E^1|$. It is long but
straightforward to show that if a finite graph $E$ has no sinks nor sources, then $L(E)\cong {\Cal C}{\Cal
K}_{A_E}(K)$.

In \cite{\AGGP, Theorem 4.1} the authors provide sufficient conditions on $A$ which yield the simplicity of ${\Cal
C}{\Cal K}_A(K)$,  in case $A$ is a finite matrix which has no row or column of zeros, and in case $A$ is not a
permutation matrix. (There is also an additional condition on an associated function $\alpha$ which must be satisfied
in order to yield the simplicity of ${\Cal C}{\Cal K}_A(K)$.)  But these conditions on $A$ eliminate both the simple
algebras $M_n(K)$ and the non-simple algebras $L(C_n)$ from consideration in \cite{\AGGP, Theorem 4.1}, since the edge
matrix for the graph given in Examples \examples\ (i) is
$$\left(\matrix 0 & 1 & 0 & \cdots & 0 \cr
               0 & 0 & 1 & \cdots & 0 \cr
          \vdots &   &   & \ddots & \vdots \cr
               0 & 0 & 0  & \cdots & 1 \cr
               0 & 0 & 0 & \cdots & 0 \cr \endmatrix\right),$$ which contains both a zero column and a zero row,
while the edge matrix for the cycle graph $C_n$ given in Example \CnExample\ is $$\left(\matrix 0 & 1 & 0 & \cdots & 0
\cr
               0 & 0 & 1 & \cdots & 0 \cr
          \vdots &   &   & \ddots & \vdots \cr
               0 & 0 & 0  & \cdots & 1 \cr
               1 & 0 & 0 & \cdots & 0 \cr \endmatrix\right),$$ which is a permutation matrix.
\endre

\vskip 10pt plus 3 pt\noindent{\bfg Acknowledgments}\par\nobreak

The authors are grateful to E. Pardo for many valuable correspondences, and to the referee for a very careful review
(especially the comments regarding the relationship between Leavitt path algebras and $C^*$-algebras). The first author
thanks P. Muhly for providing the opportunity to attend the NSF - CBMS conference on $C^*$-graph algebras held in Iowa
City, Iowa in May / June 2004. The second author was partially supported by the MCYT and Fondos FEDER,
BFM2001-1938-C02-01, MTM2004-06580-C02-02, the ``Plan Andaluz de Investigaci\'on y Desarrollo Tecnol\'ogico", FQM 336
and by a FPU fellowship by the MEC (AP2001-1368). This work was done while the second author was a Research Scholar at
the University of Colorado at Colorado Springs supported by a ``Estancias breves" FPU grant. The second author thanks
this host center for its warm hospitality.

\bigskip

%%%%%%%%%%%%%%%%%%%%%%%%%%%%%%%%%%%%%%%%%
%%%%%%%%%%%%%%%%%%%%%%%%%%  REFERENCES
%%%%%%%%%%%%%%%%%%%%%%%%%%%%%%%%%%%%%%%%%
%%%%%%%%%%%%%%%%%%%%%%%%%%%%%%%%%%%%%%%%%

\vskip 10pt plus 3 pt\noindent{\bfg References}\par\nobreak

[\AGGP] {\rms P. Ara, M.A. Gonz\'alez-Barroso, K.R. Goodearl, E. Pardo, Fractional skew monoid rings, J. Algebra 278
(2004), 104-126.}

[\BPRS] {\rms T. Bates, D. Pask, I. Raeburn, and W. Szyma\'nski, The C*-algebras of row-finite graphs, New York J.
Math. 6 (2000), 307-324.}

[\C] {\rms J. Cuntz, Simple C*-algebras generated by isometries, Comm. Math. Physics 57 (1977), 173-185.}

[\CK] {\rms J. Cuntz and W. Krieger, A class of C*-algebras and topological Markov chains, Invent. Math. 63 (1981),
25-40.}

[\KPR] {\rms A. Kumjian, D. Pask, and I. Raeburn, Cuntz-Krieger algebras of directed graphs, Pacific J. Math. 184 (1)
(1998), 161-174.}

[\Lone] {\rms W.G. Leavitt, The module type of a ring, Trans. A.M.S. 42 (1962), 113-130.}

[\Ltwo] {\rms W.G. Leavitt, The module type of homomorphic images, Duke Math. J. 32 (1965), 305-311.}

[\R] {\rms I. Raeburn, Graph algebras: operator algebras we can see, NSF - CBMS Regional Conference Series Monographs,
NSF - CBMS Conference held in Iowa City, Iowa, May 31 - June 4, 2004. (In preparation.)}

[\RS] {\rms I. Raeburn and W. Szyma\'nski, Cuntz-Krieger algebras of infinite graphs and matrices, Trans. A.M.S. 356
(1) (2004), 39-59.}

\end